\documentclass[12pt]{article}
\usepackage{amsmath,amssymb,amsthm, amsfonts}
\usepackage{etaremune, enumerate, float, verbatim}
\usepackage{algorithm}
\usepackage{algorithmic}
\usepackage{hyperref}
\usepackage{graphicx}
\usepackage{url}
\usepackage[mathlines]{lineno}
\usepackage{dsfont} 
\usepackage{tikz, graphicx, subcaption, caption}
\usepackage[algo2e,ruled,vlined]{algorithm2e}
\usepackage{mathrsfs,amsmath}
\usepackage{color}
\definecolor{red}{rgb}{1,0,0}

\definecolor{blue}{rgb}{0,0,.7}

\definecolor{green}{rgb}{0,.6,0}

\definecolor{purp}{rgb}{.5,0,.5}

\newtheorem{thm}{Theorem}[]
\newtheorem{cor}[]{Corollary}

\newtheorem{obs}[]{Observation}

\newtheorem{quest}{Question}[]

\theoremstyle{definition}

\theoremstyle{definition}
\newtheorem{defn}[]{Definition}

\theoremstyle{definition}
\newtheorem{ex}[]{Example}

\newcommand{\bit}{\begin{itemize}}
\newcommand{\eit}{\end{itemize}}
\newcommand{\ben}{\begin{enumerate}}
\newcommand{\een}{\end{enumerate}}
\newcommand{\beq}{\begin{equation}}
\newcommand{\eeq}{\end{equation}}
\newcommand{\bea}{\begin{eqnarray*}} 
\newcommand{\eea}{\end{eqnarray*}}
\newcommand{\bpf}{\begin{proof}}
\newcommand{\epf}{\end{proof}\ms}
\newcommand{\bmt}{\begin{bmatrix}}
\newcommand{\emt}{\end{bmatrix}}
\newcommand{\ms}{\medskip}

\newcommand{\noi}{\noindent}

\title{Zero forcing is not unimodal}
\author{Tony Newton
\thanks{Newton Astro Labs, London, UK (tony.newton79@gmail.com)}
\and
Boris Brimkov
\thanks{Department of Mathematics, Statistics, and Physics, Slippery Rock University, Slippery Rock, PA, USA (boris.brimkov@sru.edu)} 
}
\date{}

\begin{document}
\maketitle
\begin{abstract}
The zero forcing polynomial of a graph counts the number of distinct zero forcing sets of every possible size. In 2018, Boyer et al. asked whether the zero forcing polynomial is always unimodal, and provided some evidence that it is. We answer this question by constructing an infinite family of graphs whose zero forcing polynomials are not unimodal. Our construction also shows that the graph polynomials associated with connected forcing and total forcing are also not unimodal in general.
\end{abstract}

\noi {\bf Keywords} Zero forcing; zero forcing polynomial

\noi{\bf AMS subject classification} 05C15, 05C57, 05C31
\section{Introduction}
Given a simple graph $G$ each of whose vertices is initially colored blue or white, the \emph{zero forcing color change rule} states that a blue vertex with exactly one white neighbor causes (or \emph{forces}) that neighbor to become blue. The \emph{closure} of a set $S$ of blue vertices, denoted $\emph{cl}(S)$, is the set of blue vertices obtained after the color change rule is applied iteratively until no more color changes are possible. $S$ is a \emph{zero forcing set} of $G=(V,E)$ if $cl(S)=V$. $Z(G)$ denotes the minimum cardinality of a zero forcing set and $z(G;i)$ denotes the number of zero forcing sets of $G$ of size $i$. 

Zero forcing was independently introduced in several contexts \cite{AIM-Workshop,quantum1,fast_mixed_search} and has been widely studied. See the monograph \cite{zfbook} for a comprehensive survey of results related to zero forcing. Numerous variants of zero forcing obtained by modifying the color change rule or by adding certain restrictions to the structure of a zero forcing set have also been studied. For example, \emph{connected forcing} requires the zero forcing set to induce a connected subgraph~\cite{cf}, and \emph{total forcing} requires the zero forcing set to induce a graph with no isolated vertices~\cite{tf}. 

The \emph{zero forcing polynomial} of a graph $G$ on $n$ vertices is defined as
\begin{equation}
\label{e1}
\mathcal{Z}(G;x)=\sum_{i=1}^n z(G;i) x^i.
\end{equation}
\noindent In \cite{BB}, the authors asked whether the zero forcing polynomial of every graph is unimodal (i.e., whether there is always some $i^*\in\{1,\ldots,n\}$ such that $z(G;1)\leq\ldots\leq z(G;i^*)\geq \ldots\geq z(G;n)$). They showed that the first half of the coefficients are always non-decreasing, i.e., that  $z(G;i)\leq z(G;i+1)$ for $1\leq i< \frac{n}{2}$. They also showed that the zero forcing polynomials of several families of graphs including paths, cycles, complete graphs, and complete multipartite graphs are unimodal. Thus, it seemed plausible that the answer to the question is affirmative. In this note, we settle this question by constructing an infinite family of graphs whose zero forcing polynomials are not unimodal. Our construction also shows that the graph polynomials associated with connected forcing and total forcing (respectively denoted $\mathcal{Z}_c(G;x)$ and $\mathcal{Z}_t(G;x)$, and defined analogously to \eqref{e1}) are also not unimodal in general.

\subsection*{Preliminaries}
A simple graph $G=(V,E)$ consists of a vertex set $V$ and an edge set $E$ of two-element subsets of $V$.  The \emph{order} of $G$ is denoted by $n=|V|$. Two vertices $v,w\in V$ are \emph{adjacent}, or \emph{neighbors}, if $\{v,w\}\in E$. The \emph{neighborhood} of $v\in V$, denoted $N(v)$, is the set of all vertices which are adjacent to $v$. The \emph{degree} of $v\in V$ is $\deg(v)=|N(v)|$. An \emph{isolated vertex} is a vertex that is not adjacent to any other vertex, and a \emph{universal vertex} is a vertex that is adjacent to all other vertices. 
Given $S \subset V$, the \emph{induced subgraph} $G[S]$ is the subgraph of $G$ whose vertex set is $S$ and whose edge set consists of all edges of $G$ which have both endpoints in $S$. 

The complete graph on $n$ vertices is denoted by $K_n$ and the empty graph on $n$ vertices is denoted by $\overline{K}_n$. 
The \emph{join} of disjoint graphs $G_1$ and $G_2$, denoted $G_1\vee G_2$, is the graph obtained by adding an edge between every vertex of $G_1$ and every vertex of $G_2$. For other graph theoretic terminology and definitions, we refer the reader to~\cite{west}.

Two distinct vertices $u$ and $v$ in $G$ are \emph{twins} if
$N(u)\setminus\{v\}=N(v)\setminus\{u\}$, i.e., $u$ and $v$ have
identical neighborhoods outside $\{u,v\}$. A \emph{fort} of a graph $G=(V,E)$ is a non-empty set $F\subset V$ such that no vertex outside $F$ is adjacent to exactly one vertex in $F$. It was shown in \cite{brimkov_zf1} that every zero forcing set of a graph intersects every fort of the graph. Since each pair of twin vertices is a fort, it follows that every zero forcing set intersects every pair of twin vertices. This can be generalized to larger sets of pairwise twins as follows. 
\begin{obs}
\label{ob1}
Let $G=(V,E)$ be a graph and $S\subset V$ be a set such that every pair of vertices in $S$ are twins. Then, every zero forcing set of $G$ contains all or all-but-one vertices of $S$.   
\end{obs}

\section{Main result}
\begin{figure}[h!]
    \centering
    \resizebox{0.7\textwidth}{!}{
    \begin{tikzpicture}[
        vertex/.style={circle, draw=black, thick, minimum size=0.65cm, inner sep=0pt},
        line/.style={thick}
    ]
    \node[vertex] (A) at (0,2) {$a_1$};
    \node[vertex] (B) at (0,0) {$a_2$};
    \node[vertex] (C) at (0,-2) {$a_3$};
    \node[vertex] (D) at (4,1) {$b_1$};
    \node[vertex] (E) at (4,-1) {$b_2$};
    \node[vertex] (F) at (8,1) {$c_1$};
    \node[vertex] (G) at (8,-1) {$c_2$};
    \node[vertex] (H) at (12,2) {$d_1$};
    \node[vertex] (I) at (12,0) {$d_2$};
    \node[vertex] (J) at (12,-2) {$d_3$};
    \draw[line] (A)--(D);
    \draw[line] (A)--(E);
    \draw[line] (B)--(D);
    \draw[line] (B)--(E);
    \draw[line] (C)--(D);
    \draw[line] (C)--(E);
    \draw[line] (D)--(E);
    \draw[line] (D)--(F);
    \draw[line] (D)--(G);
    \draw[line] (E)--(F);
    \draw[line] (E)--(G);
    \draw[line] (F)--(H);
    \draw[line] (F)--(I);
    \draw[line] (F)--(J);
    \draw[line] (G)--(H);
    \draw[line] (G)--(I);
    \draw[line] (G)--(J);

    \end{tikzpicture}
    }
    \caption{The graph $H$.}
    \label{fig1}
\end{figure}
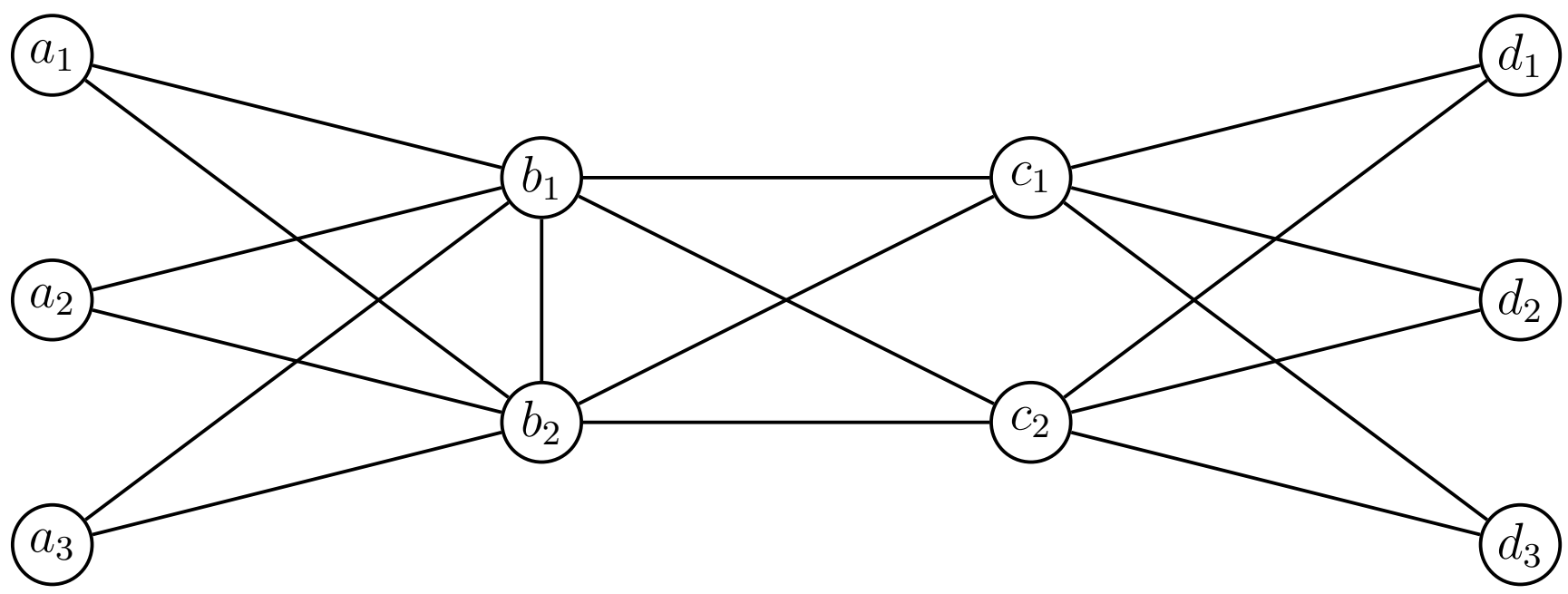

\begin{defn}
\label{d1}
Let $H$ be the graph shown in Figure \ref{fig1}, with vertices labeled as in the figure. Let $A=\{a_1,a_2,a_3\}$,
$B=\{b_1,b_2\}$,
$C=\{c_1,c_2\}$,
$D=\{d_1,d_2,d_3\}$ be a partition of the vertices of $H$. 
\end{defn}

\begin{defn}
\label{d2}
Let $G_{p,q,r}=H\vee \overline K_p \vee \overline K_q \vee K_r$. Let $P$, $Q$, and $R$ respectively denote the vertex sets of the copies of $\overline K_p$, $\overline K_q$ and $K_r$ used in the construction of $G_{p,q,r}$, and let $H'=P\cup Q\cup R$. Let $G_t=G_{t,t,t}$ and let $n=10+3t$ denote the order of $G_t$.
\end{defn}

\begin{figure}[!ht]
    \centering
\includegraphics[width=0.8\linewidth]{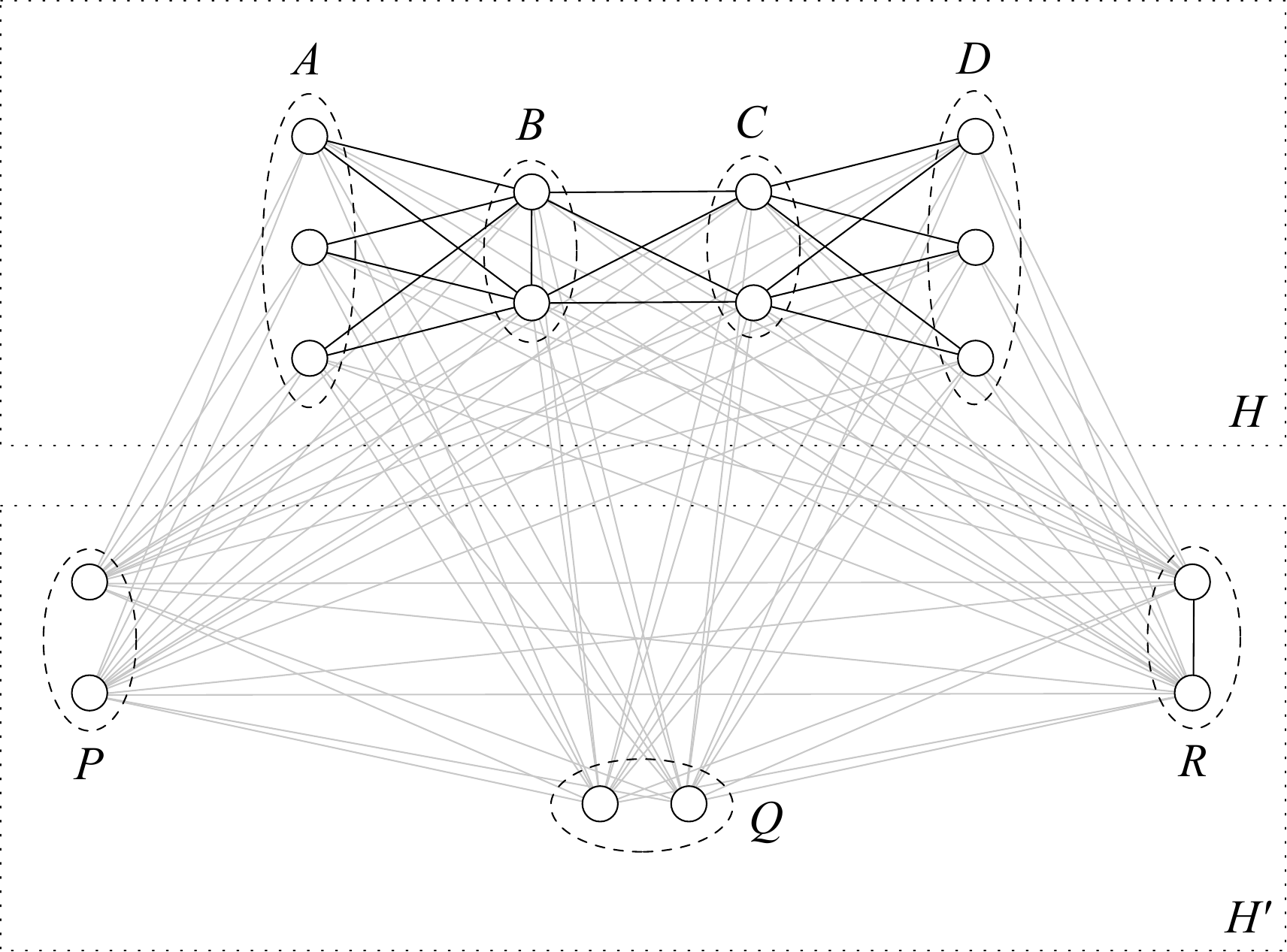}
    \caption{The graph $G_2$.}
    \label{fig2}
\end{figure}

\begin{ex}
The graph $G_2$ is shown in Figure \ref{fig2}.
\end{ex}

\begin{thm}
\label{main_thm}
For $t\geq 2$, $\mathcal{Z}(G_t;x)=x^n+(10+3t)x^{n-1}+(37+30t+3t^2)x^{n-2}+(60+99t)x^{n-3}+(36+108t)x^{n-4}$.
\end{thm}
\proof
Let $S$ be an arbitrary zero forcing set of $G_t$. Let $X$ be any one of the sets $A$, $B$, $C$, $D$, $P$, $Q$, $R$ that partition the vertices of $G_t$ according to Definitions \ref{d1} and \ref{d2}. By construction, every pair of vertices $x_1,x_2\in X$ are twins. Thus, by Observation \ref{ob1}, $S$ contains either all or all-but-one vertices of $X$. Let $s=|H'\setminus S|$ be the number of vertices of $H'$ that are not in $S$. Since $S$ can exclude at most one vertex from each of $P$, $Q$, and $R$, it follows that $s\in \{3,2,1,0\}$. We examine each of these cases separately. 

\emph{Case 1: $s=3$.} This means $S$ excludes exactly one vertex from each of $P$, $Q$, and $R$. However, this implies that no vertex can perform the first force, since every vertex in $G_t$ would be adjacent to at least two white vertices. Thus, this case cannot occur.

\emph{Case 2: $s=2$.} If $S$ excludes any vertex of $H$ in addition to the two vertices in $H'$, then again every vertex in $G_t$ would be adjacent to at least two white vertices and no vertex could perform the first force.
Thus, if $s=2$, it follows that $S$ must include every vertex of $H$. On the other hand, any set $S$ that excludes exactly two vertices $u,v\in H'$ (each of which belongs to a different set among $P$, $Q$, $R$) and includes every other vertex of $G_t$ is a zero forcing set. To see why, note that by Pigeonhole Principle, $(P\cup Q)\cap \{u,v\}\neq \emptyset$. Without loss of generality, suppose $u\in P$ and let $w$ be another vertex in $P$. Then, $w$ can force $v$ and $v$ can force $u$.

\emph{Case 3: $s=1$.} Suppose for contradiction that in addition to the vertex in $H'$, $S$ also excludes a vertex in $B$ and a vertex in $C$. Then, every vertex in $G_t$ would be adjacent to at least two white vertices and no vertex could perform the first force. Thus, if $s=1$, $S$ cannot exclude a vertex from both $B$ and $C$. On the other hand, any set $S$ that excludes exactly one vertex $u\in H'$, excludes at most one vertex from each of $A$, $B$, $C$, and $D$, and does not simultaneously exclude a vertex from both $B$ and $C$, is a zero forcing set. To see why, suppose without loss of generality that $S$ either excludes some (possibly non-proper) subset of $\{u, a_1,b_1,d_1\}$, or of $\{u, a_1,c_1,d_1\}$. In the first case, $d_2$ can force $u$ in the first timestep, then $a_2$ can force $b_1$, then $b_1$ can force $a_1$, and then $c_1$ can force $d_1$ (if $a_1$ and/or $b_1$ and/or $d_1$ do not already belong to $S$). In the second case, $a_2$ can force $u$ in the first timestep, then $c_2$ can force $d_1$, then $d_1$ can force $c_1$, and then $b_1$ can force $a_1$ (if $a_1$ and/or $c_1$ and/or $d_1$ do not already belong to $S$).

\emph{Case 4: $s=0$.} This means $S$ includes every vertex in $H'$ and  excludes at most one vertex from each of $A$, $B$, $C$, and $D$. Any such set is a zero forcing set. To see why, suppose without loss of generality that $S$ excludes some (possibly non-proper) subset of the vertices $\{a_1, b_1, c_1, d_1\}$. In the first timestep, $a_2$ can force $b_1$ and $d_2$ can force $c_1$ (if $b_1$ and/or $c_1$ do not already belong to $S$). Then, $b_1$ can force $a_1$ and $c_1$ can force $d_1$ (if $a_1$ and/or $d_1$ do not already belong to $S$).

Since the cases above hold for any zero forcing set, they also hold for a minimum zero forcing set. Thus, a minimum zero forcing set can exclude at most 4 vertices of $G_t$, and so $Z(G_t)=n-4$.

We will now count $z(G_t;i)$ for each $i\in \{n-4,\ldots,n\}$. Clearly, $z(G_t;n)=1$. Since $G_t$ is connected, it is also clear that $z(G_t;n-1)=n=10+3t$. 
To count $z(G_t;n-2)$, note that there are three possibilities:
\begin{itemize}
\item Both of the excluded vertices are in $H$. By Case 4, the excluded vertices can be in any of the pairs of sets $(A,B), (A,C), (A,D), (B,C), (B,D), (C,D)$, which respectively yields $3\cdot 2+3\cdot 2+3\cdot 3+2\cdot 2+2\cdot 3+2\cdot 3=37$ distinct zero forcing sets. 

\item One of the excluded vertices is in $H$ and the other is in $H'$. Since there is only a single excluded vertex in $H$, by Case 3 the excluded vertices can be in any of the pairs of sets $(A,H'), (B,H'), (C,H'), (D,H')$, which respectively yields $3\cdot 3t+2\cdot 3t+2\cdot 3t+3\cdot 3t=30t$ distinct zero forcing sets.  

\item Both of the excluded vertices are in $H'$. Since there are no excluded vertices in $H$, by Case 2 the excluded vertices can be in any of the pairs of sets $(P,Q), (P,R), (Q,R)$, which respectively yields $t\cdot t+t\cdot t+t\cdot t=3t^2$ distinct zero forcing sets. 
\end{itemize}
Thus, in total $z(G_t;n-2)=37+30t+3t^2$. Similarly, to count $z(G_t;n-3)$,  there are two possibilities:

\begin{itemize}
\item All three of the excluded vertices are in $H$. By Case 4, the excluded vertices can be in any of the triples of sets $(A,B,C), (A,B,D), (A,C,D), (B,C,D)$, which respectively yields $3\cdot 2\cdot 2+3\cdot 2\cdot 3+3\cdot 2\cdot 3+2\cdot 2\cdot 3=60$ distinct zero forcing sets. 

\item Two of the excluded vertices are in $H$ and the third is in $H'$. By Case 3, the excluded vertices can be in any of the triples of sets $(A,B,H')$, $(A,C,H')$, $(A,D,H')$, $(B,D,H')$, $(C,D,H')$, which respectively yields $3\cdot 2\cdot 3t+3\cdot 2\cdot 3t+3\cdot 3\cdot 3t+2\cdot 3\cdot 3t+2\cdot 3\cdot 3t=99t$ distinct zero forcing sets.  
\end{itemize}
Note that it is not possible for two or three of the excluded vertices to be in $H'$ by Cases 2 and 1, respectively. Thus, in total $z(G_t;n-3)=60+99t$. Finally, to count $z(G_t;n-4)$, there are again two possibilities:
\begin{itemize}
\item All four of the excluded vertices are in $H$, with one in each of the sets $A,B,C,D$. This yields $3\cdot 2\cdot 2\cdot 3=36$ distinct zero forcing sets. 

\item Three of the excluded vertices are in $H$ and the fourth is in $H'$. By Case 3, the excluded vertices can be in any of the triples of sets $(A,B,D,H'), (A,C,D,H')$, which respectively yields $3\cdot 2\cdot 3\cdot 3t+3\cdot 2\cdot 3\cdot 3t=108t$ distinct zero forcing sets.  
\end{itemize}
Again, by Cases 2 and 1, it is not possible for more than one of the excluded vertices to be in $H'$. Thus, in total $z(G_t;n-4)=36+108t$. Plugging in the above coefficients into \eqref{e1} , we obtain $\mathcal{Z}(G_t;x)=x^n+(10+3t)x^{n-1}+(37+30t+3t^2)x^{n-2}+(60+99t)x^{n-3}+(36+108t)x^{n-4}$.
\qed
\begin{ex}
$\mathcal{Z}(G_2;x)=x^{16}+16x^{15}+109x^{14}+258x^{13}+ 252x^{12}.$
\end{ex}

\begin{cor}
\label{c1}
For $t\geq 24$, $\mathcal{Z}(G_t;x)$ is not unimodal.
\end{cor}
\proof
We will show that for $t\geq 24$, $z(G_t,n-4)>z(G_t,n-3)$ and $z(G_t,n-3)<z(G_t,n-2)$. Indeed, using the coefficients derived in Theorem \ref{main_thm} we have 
\begin{align*}
 z(G_t,n-4)-z(G_t,n-3)&=\\(36+108t)-(60+99t)&=\\9t-24&>0
\end{align*}
and
\begin{align*} z(G_t,n-2)-z(G_t,n-3)&=\\
 (37+30t+3t^2)-(60+99t)
 &=\\
 3(t-24)^2+75(t-24)+49&>0.\qedhere
\end{align*}
\begin{ex}
$G_{24}$ is a graph on 82 vertices with zero forcing polynomial $\mathcal{Z}(G_{24};x)=x^{82}+82x^{81}+2485x^{80}+2436x^{79}+ 2628x^{78}$, which is not unimodal.
\end{ex}

\begin{cor}
For $t\geq 24$, $\mathcal{Z}_c(G_t;x)$ and $\mathcal{Z}_t(G_t;x)$ are not unimodal.
\end{cor}
\proof
As shown in the proof of Theorem \ref{main_thm}, every zero forcing set of $G_t$ contains either all or all-but-one vertices of $R$. Since $|R|=t\geq 24$ it follows that every zero forcing set $S$ of $G_t$ includes a vertex of $R$ which is a universal vertex, and hence that $S$ induces a connected subgraph with no isolated vertices. 
On the other hand, every connected or total forcing set of a graph $G$ must also be a zero forcing set of $G$.
Thus, $S$ is a zero forcing set of $G_t$ if and only if $S$ is a connected forcing set of $G_t$ if and only if $S$ is a total forcing set of $G_t$. It follows that 
$\mathcal{Z}(G_t;x)=\mathcal{Z}_c(G_t;x)=\mathcal{Z}_t(G_t;x)$, and so by Corollary \ref{c1}, $\mathcal{Z}_c(G_t;x)$ and $\mathcal{Z}_t(G_t;x)$ are not unimodal for $t\geq 24$.
\qed

\section{Conclusion}
In this note we showed that the zero forcing polynomial of a graph is in general not unimodal. It would be interesting to explore whether certain graph properties guarantee unimodality or non-unimodality of the zero forcing polynomial. Specifically, we pose the following open question.
\begin{quest}
Is the zero forcing polynomial of a tree always unimodal?
\end{quest}
\noindent There are smaller graphs than $G_{24}$ whose zero forcing polynomials are not unimodal. For example, it can be shown with a similar argument as Theorem \ref{main_thm} that $\mathcal{Z}(G_{24,24,23};x)=2592x^{77}+2403x^{78}+2427x^{79}+81x^{80}+x^{81}$. The following remains an open question.
\begin{quest}
What is the smallest graph with a non-unimodal zero forcing polynomial?
\end{quest} 
\noindent Finally, it would be interesting to determine whether non-unimodality of the zero forcing polynomial is rare (e.g., among all graphs of order $n$), or whether there are in fact many more examples to be found. Some work in this direction has been initiated in \cite{almost_uni}.
\begin{quest}
Do almost all graphs have unimodal zero forcing polynomials?
\end{quest}

\section*{AI Disclosure}
Generative AI systems were used during the mathematical investigation and to generate an early draft of the manuscript. The final draft, including all proofs, was rewritten from scratch entirely by humans.

\end{document}